\documentclass[twoside,reqno,A4]{amsart}

\usepackage{latexsym}

\newcommand{\xquad}{\quad\qquad\quad}
\newcommand{\qdn}{\hspace*{-1.5mm}}
\newcommand{\qqdn}{\hspace*{-2.5mm}}

\newcommand{\dst}{\displaystyle}
\newcommand{\sst}{\scriptstyle}

\newcommand{\+}{&\qqdn}%

\newcommand{\ang}[1]{\langle{#1}\rangle}

\newcommand{\mb}[1]{\mathbb{#1}}
\newcommand{\mc}[1]{\mathcal{#1}}

\newcommand{\ffnk}[4]{\left[\qdn\ba{#1}#3\\#4\ea{\!\Big|\:#2}\right]}

\newcommand{\binm}{\binom}
\newcommand{\sbnm}[2]{\Bigl(\!\ba{c}\!#1\!\\#2\ea\!\Bigr)}

\newcommand{\double}[2]{\genfrac{}{}{0mm}{1}{#1}{#2}}

\newcommand{\be}{\begin{equation}}
\newcommand{\ee}{\end{equation}}
\newcommand{\ba}{\begin{array}}
\newcommand{\ea}{\end{array}}
\newcommand{\bmn}{\begin{eqnarray}}
\newcommand{\emn}{\end{eqnarray}}
\newcommand{\bnm}{\begin{eqnarray*}}
\newcommand{\enm}{\end{eqnarray*}}
\newcommand{\bln}{\begin{subequations}}
\newcommand{\eln}{\end{subequations}}

\newcommand{\eqn}[1]{\begin{equation}#1\end{equation}}
\newcommand{\xalignx}[1]{\begin{align}#1
            \end{align}}     
\newcommand{\xalignz}[1]{\begin{align*}#1
            \end{align*}}    

\newcommand{\xalignedx}[1]{\begin{aligned}#1
            \end{aligned}}  
\newcommand{\mult}[2]{\begin{array}{#1}#2\end{array}}%

\newcommand{\centro}[1]
           {\begin{center}#1\end{center}}

\newcommand{\lam}{\lambda}

\newtheorem{thm}{Theorem}

\newcommand{\thank}[2]{\begin{center}\parbox{#1}
{{\sc\bf Acknowledgement}:\,{\small\it#2}}\end{center}}

\newcommand{\bibtm}[4]{\bibitem{kn:#1}{#2,}~\emph{#3,}~{#4.}}	
\newcommand{\cito}[1]{\cite{kn:#1}}	
\newcommand{\citu}[2]{\cite[#2]{kn:#1}}

\begin{document} 

\title{Hypergeometric Series and \\ Harmonic Number Identities}
\author{CHU Wenchang - DE DONNO Livia}

\address{
 	Dipartimento\: \:di\: \:Matematica	\newline
     	Universit\a`{a} degli Studi di Lecce\newline
     	Lecce-Arnesano\: \:P.\:O.\:Box\:193	\newline
      73100 Lecce,$\quad$ITALIA		\newline              	
	tel 39+0832+297409			\newline 
	fax 39+0832+297594			\newline              	
	Email \emph{chu.wenchang@unile.it}	\newline
	Email \emph{liviadedonno@libero.it}}

\begin{abstract}
The classical hypergeometric summation theorems
are exploited to derive several striking identities
on harmonic numbers including those discovered recently 
by Paule and Schneider (2003).
\end{abstract} 

\maketitle
\markboth{Chu - De Donno: Hypergeometric Series 
		and Harmonic Number Identities}
         {Chu - De Donno: Hypergeometric Series 
		and Harmonic Number Identities}
\thispagestyle{empty}


\vspace*{-5mm}\section{Introduction and Notation}

Let $x$ be an indeterminate. The generalized harmonic numbers 
are defined to be partial sums of the harmonic series: 
\eqn{H_0(x)=0 
\quad\text{and}
\quad H_ n(x)
\:=\:\sum_{k=1}^n 
\frac{1}{x+k}
\quad\text{for}\quad 
n=1,2,\cdots.} 
For $x=0$ in particular, they reduce to the classical 
harmonic numbers:
\eqn{H_0=0
\quad\text{and}\quad 
H_n=\sum_{k=1}^n\frac{1}{k}
\quad\text{for}\quad 
n=1,2,\cdots.}

Given a differentiable function $f(x)$, 
denote two derivative operators by
\[\mc{D}_xf(x)\:=\:\frac{d}{dx}f(x)
\quad\text{and}\quad
\mc{D}_0f(x)\:=\:\frac{d}{dx}f(x)\Big|_{x=0}.\]

Then it is an easy exercise to compute 
the derivative of binomial coefficients  
\[\mc{D}_x\sbnm{x+n}{m}
\:=\:\sbnm{x+n}{m}
\sum_{\ell=1}^m
\frac{1}{1+x+n-\ell}\]
which can be stated in terms of the generalized harmonic numbers as
\eqn{\mc{D}_x\sbnm{x+n}{m}
\:=\:\sbnm{x+n}{m}
\big\{H_n(x)-H_{n-m}(x)\big\}, 
\quad (m\le n).}
In this paper, we will frequently use its evaluation at $x=0$:
\eqn{\mc{D}_0\binm{x+n}{m}
\:=\:\label{d+bin}
\binm{n}{m}
\big\{H_n-H_{n-m}\big\},
\qquad (m\le n).}

For the inverse binomial coefficients,  
the analogous results read as
\[\mc{D}_x\sbnm{x+n}{m}^{-1}
\:=\:\sbnm{x+n}{m}^{-1}
\sum_{\ell=1}^m
\frac{-1}{1+x+n-\ell}\]
and the explicit harmonic number expressions
\xalignx{
\mc{D}_x\binm{x+n}{m}^{-1}
&\:=\:\sbnm{x+n}{m}^{-1}
\big\{H_{n-m}(x)-H_n(x)\big\}, 
&&(m\le n)\\
\mc{D}_0\binm{x+n}{m}^{-1}
&\:=\:\label{d-bin}
\binm{n}{m}^{-1}
\big\{H_{n-m}-H_n\big\},
&\quad&(m\le n).}

As pointed out by Richard Askey (cf.\:\cito{andrews} and \cito{paule}), 
expressing harmonic numbers in terms of differentiation of binomial 
coefficients can be traced back to Issac Newton. Following the work 
of the two papers cited above, we will explore further the application 
of derivative operators to hypergeometric summation formulas. Several 
striking harmonic number identities discovered in \cito{paule} will 
be recovered and some new ones will be established. 

Because hypergeometric series will play a central role in the present work,
we reproduce its notation for those who are not familiar with it.
Roughly speaking, a hypergeometric series is a series $\sum C_n$
where the term ratio $C_{n+1}/C_n$ is a rational function in $n$. 
If the shifted factorial is defined by  
\eqn{(c)_0=1
\quad\text{and}\quad
(c)_n=c(c+1)\cdots(c+n-1)
\quad\text{for}\quad
n=1,2,\cdots}
then the hypergeometric series (cf.\:\cito{bailey}) 
reads explicitly as
\eqn{_{1+p}F_{q}
\ffnk{cccc}{z}{a_0,&a_1,\cdots,&a_p}
     {&b_1,\cdots,&b_q}=
\sum_{n=0}^{\infty}
\frac{(a_0)_n(a_1)_n\cdots(a_p)_n}
     {n!\:(b_1)_n\cdots(b_q)_n}\:z^n.}
In order to illustrate how to discover harmonic number 
identities from hypergeometric series, we start with 
the Chu-Vandermonde-Gauss formula~\citu{bailey}{\S1.3}: 
\[_{2}F_{1}\ffnk{cc}{1}{-n,&a}{&c}
=\frac{(c-a)_n}{(c)_n}.\]
Under parameter replacements
$a\to-n-\mu n$
and $c\to 1+\lam n+x$
with $\lam,\:\mu\in\mb{N}_0$,
it can equivalently be stated as the following 
binomial convolution identity
\eqn{\sum_{k=0}^n
\binm{n+\mu n}{k}
\binm{x+\lam n+n}{n-k}
\:=\:\label{eq:chu}
\binm{x+\lam n+\mu n+2n}{n}.}
In view of \eqref{d+bin}, we derive, by applying 
the derivative operator $\mc{D}_0$ to both sides 
of the last identity, the following relation:   
\[\sum_{k=0}^n
\sbnm{n\!+\!\mu n}{k}
\sbnm{n\!+\!\lam n}{n\!-\!k}
\Big\{H_{\lam n+n}-H_{\lam n+k}\Big\}=
\sbnm{2n\!+\!\lam n\!+\!\mu n}{n}
\big\{\!H_{\lam n+\mu n+2n}
\!-\!H_{\lam n+\mu n+n}\!\big\}.\]
According to the factor inside the braces $\{\cdots\}$,
spliting the left hand side into two sums with respect 
to $k$ and then evaluating the first one by \eqref{eq:chu},  
we get immediately the following simplified result.

\begin{thm}\label{chu}
With $\lam,\:\mu\in\mb{N}_0$, there holds
the following harmonic number identity:
\[\sum_{k=0}^n
\sbnm{n\!+\!\mu n}{k}
\sbnm{n\!+\!\lam n}{n\!-\!k}
H_{\lam n+k}=
\sbnm{2n\!+\!\lam n\!+\!\mu n}{n}
\big\{\!H_{\lam n+n}\!+\!H_{\lam n+\mu n+n}
\!-\!H_{\lam n+\mu n+2n}\!\big\}.\]
\end{thm}
One interesting special case corresponding to $\mu=0$
can be stated as 
\eqn{\sum_{k=0}^n
\sbnm{n}{k}\sbnm{n+\lam n}{n-k}
H_{\lam n+k}
\:=\:\label{wench}
\sbnm{2n+\lam n}{n}
\big\{2H_{\lam n+n}
-H_{\lam n+2n}\big\}.}
It can be further specialized, with $\lam=0$, to  
\eqn{\sum_{k=0}^n\sbnm{n}{k}^2H_{k}
\:=\:\sbnm{2n}{n}
\big\{2H_{n}-H_{2n}\big\}.}
There exist numerous hypergeometric series identities.
However we are not going to have a full coverage about
how they can be used to find harmonic number identities. 
The authors will limit themselves to examine, by the derivative 
operator method, only the classical identities named 
after Pfaff-Saalsch\"utz, Dougall-Dixon and the Whipple 
transformation in next three sections. As applications,
we will tabulate $26$ closed formulas and $21$ transformations 
on harmonic numbers at the end of the paper. 

Just like the demonstration of Theorem~\ref{chu} and \eqref{wench},
we will examine the above-mentioned hypergeometric theorems 
in the three steps: reformulation in terms of binomial formulas,  
application of the derivative operator $\mc{D}_0$ and reduction
to harmonic number identities by specifying parameters.
Because all the computations involved in the paper are routine 
manipulations on finite series, we will therefore omit the details 
for the limit of space.


\section{The Pfaff-Saalsch{\"u}tz Theorem}

Recall the Saalsch{\"u}tz theorem~\citu{bailey}{\S 2.2} 
\[{_{3}F_{2}
\ffnk{ccc}{1}{-n,\+a,\+b}
     {\+c,\+1+a+b-c-n}
\:=\:\frac{(c-a)_n(c-b)_n}
          {(c)_n(c-a-b)_n}}.\]
Performing the parameter replacement 
\[\left.\xalignedx{
a&\to-n-\mu n-\mu'x\\
b&\to1+\lam n+\lam'x\\
c&\to1+\nu n+\nu'x}
\right\}\qquad
(\lam,\mu,\nu\in\mb{N}_0)\]
we may express it as a binomial identity
\[\sum_{k=0}^n
\frac{\binm{n}{k}
	\binm{k+\lam n+\lam'x}{k}
	\binm{n+\mu n+\mu'x}{k}}
     {\binm{k+\nu n+\nu'x}{k}
	\binm{k+(\lam-\mu-\nu-2)n+(\lam'-\mu'-\nu')x}{k}}
\:=\:
\frac{\binm{(\lam-\nu)n+(\lam'-\nu')x}{n}
	\binm{(\mu+\nu+2)n+(\mu'+\nu')x}{n}}
     {\binm{n+\nu n+\nu'x}{n}
	\binm{(\lam-\mu-\nu-1)n+(\lam'-\mu'-\nu')x}{n}}.\]
Applying $\mc{D}_0$ to the cases $\mu'=\nu'=0$, $\lam'=\nu'=0$
and $\lam'=\mu'=0$ of the last identity, we get respectively
the following harmonic number identities.
\begin{thm}
For $\lam,\:\mu,\:\nu\in\mb{N}_0$ with $\lam>1+\mu+\nu$, 
we have the harmonic number identity:\label{thm:ps-lam}
\xalignz{&\sum_{k=0}^n
\frac{\binm{n}{k}
	\binm{\lam n+k}{k}
	\binm{\mu n+n}{k}}
     {\binm{\nu n+k}{k}
	\binm{(\lam-\mu-\nu-2)n+k}{k}}
\big\{H_{\lam n+k}
-H_{(\lam-\mu-\nu-2)n+k}\big\}\\
&\:=
\frac{\binm{(\lam-\nu)n}{n}
	\binm{(\mu+\nu+2)n}{n}}
     {\binm{\nu n+n}{n}
	\binm{\lam-\mu-\nu-1)n}{n}}
\Big\{\qdn\mult{c}
{H_{(\lam-\nu)n}-H_{(\lam-\nu-1)n}\\
+H_{\lam n}-H_{(\lam-\mu-\nu-1)n}}
\qdn\Big\}.}
\end{thm}

\begin{thm}
For $\lam,\:\mu,\:\nu\in\mb{N}_0$ with $\lam>1+\mu+\nu$, 
we have the harmonic number identity:\label{thm:ps-mu}
\xalignz{
&\sum_{k=0}^n
\frac{\binm{n}{k}
	\binm{\lam n+k}{k}
	\binm{\mu n+n}{k}}
     {\binm{\nu n+k}{k}
	\binm{(\lam-\mu-\nu-2)n+k}{k}}
\big\{H_{\mu n\!+\!n\!-\!k}
-H_{(\lam\!-\!\mu\!-\!\nu\!-\!2)n+k}\big\}\\
&\:=
\frac{\binm{(\lam-\nu)n}{n}
	\binm{(\mu+\nu+2)n}{n}}
     {\binm{\nu n+n}{n}
	\binm{\lam-\mu-\nu-1)n}{n}}
\Big\{\qdn\mult{c}
{H_{(\mu+\nu+1)n}-H_{(\mu+\nu+2)n}\\
+H_{\mu n+n}-H_{(\lam-\mu-\nu-1)n}
}\qdn\Big\}.}
\end{thm}

\begin{thm}
For $\lam,\:\mu,\:\nu\in\mb{N}_0$ with $\lam>1+\mu+\nu$, 
we have the harmonic number identity:\label{thm:ps-nu}
\xalignz{
&\sum_{k=0}^n
\frac{\binm{n}{k}
	\binm{\lam n+k}{k}
	\binm{\mu n+n}{k}}
     {\binm{\nu n+k}{k}
	\binm{(\lam-\mu-\nu-2)n+k}{k}}
\big\{H_{\nu n+k}
-H_{(\lam-\mu-\nu-2)n+k}\big\}\\
&\:=
\frac{\binm{(\lam-\nu)n}{n}
	\binm{(\mu+\nu+2)n}{n}}
     {\binm{\nu n+n}{n}
	\binm{\lam-\mu-\nu-1)n}{n}}
\Bigg\{\qdn\mult{c}
{H_{(\mu+\nu+1)n}-H_{(\mu+\nu+2)n}\\
+H_{(\lam-\nu)n}-H_{(\lam-\nu-1)n}\\
+H_{\nu n+n}-H_{(\lam-\mu-\nu-1)n}
}\qdn\Bigg\}.}
\end{thm}


\section{The Dougall-Dixon Theorem}
This section will explore the Dougall-Dixon 
theorem~\citu{bailey}{\S 4.3}
\[_{5}F_{4}
\ffnk{ccccc}{1}{a,\+1+a/2,\+b,\+d,\+-n}
     {\+a/2,\+1+a-b,\+1+a-d,\+1+a+n}
=\frac{(1+a)_n(1+a-b-d)_n}
      {(1+a-b)_n(1+a-d)_n}\]
to establish harmonic number identities.


\subsection{}
Performing parameter replacement   
\[\left.\xalignedx{
a&\to-n-x\\
b&\to1+bn\\
d&\to1+dn}
\right\}\quad
(b,\:d\in\mb{N}_0)\]
we can reformulate the Dougall-Dixon theorem
as the following binomial identity: 
\[\sum_{k=0}^n
\big\{x\!+\!n\!-\!2k\big\}
\sbnm{n}{k}
\frac{\binm{x+n}{k}
\binm{k+bn}{k}
\binm{k+dn}{k}}
{\binm{k-x}{k}
\binm{x+bn+n}{k}
\binm{x+dn+n}{k}}
=x
\frac{\binm{x+n}{n}
\binm{1+x+bn+dn+n}{n}}
{\binm{x+bn+n}{n}
\binm{x+dn+n}{n}}\]
which leads us, under the derivative operator 
$\mc{D}_0$, to the following result.
\begin{thm}\label{thm:2/0} With $b,\:d\in\mb{N}_0$, 
there holds the following harmonic number identity:
\[\sum_{k=0}^n
\sbnm{n}{k}^2
\frac{\binm{k+bn}{k}
\binm{k+dn}{k}}
{\binm{n+bn}{k}
\binm{n+dn}{k}}
\Big\{1+(n\!-\!2k)
\big(2H_k\!-\!H_{bn+k}
\!-\!H_{dn+k}\big)\Big\}
=\frac{\binm{1+bn+dn+n}{n}}
      {\binm{n+bn}{n}
       \binm{n+dn}{n}}.\]
\end{thm}


\subsection{}
Performing parameter replacement   
\[\left.\xalignedx{
a&\to-n-x\\
b&\to1+bn\\
d&\to-n-dn}
\right\}\quad
(b,\:d\in\mb{N}_0)\]
we can reformulate the Dougall-Dixon theorem
as the following binomial identity: 
\[\sum_{k=0}^n
\big\{x\!+\!n\!-\!2k\big\}
\sbnm{n}{k}
\frac{\binm{x+n}{k}
\binm{k+bn}{k}
\binm{n+dn}{k}}
{\binm{k-x}{k}
\binm{n+bn+x}{k}
\binm{k+dn-x}{k}}
=(-1)^nx
\frac{\binm{x+n}{n}
\binm{x+bn-dn}{n}}
{\binm{n+bn+x}{n}
\binm{n+dn-x}{n}}\]
which leads us, under the derivative operator 
$\mc{D}_0$, to the following result.
\begin{thm}\label{thm:1/1} With $b,\:d\in\mb{N}_0$, 
there holds the following harmonic number identity:
\[\sum_{k=0}^n
\sbnm{n}{k}^2
\frac{\binm{k+bn}{k}
\binm{n+dn}{k}}
{\binm{n+bn}{k}
\binm{k+dn}{k}}
\Big\{1+(n\!-\!2k)
\big(2H_k\!-\!H_{bn+k}
\!+\!H_{dn+k}\big)\Big\}
=(-1)^n
\frac{\binm{bn-dn}{n}}
{\binm{n+bn}{n}
\binm{n+dn}{n}}.\]
\end{thm}


\subsection{}
Performing parameter replacement   
\[\left.\xalignedx{
a&\to-n-x\\
b&\to-n-bn\\
d&\to-n-dn}
\right\}\quad
(b,\:d\in\mb{N}_0)\]
we can reformulate the Dougall-Dixon theorem
as the following binomial identity: 
\[\sum_{k=0}^n
\big\{x\!+\!n\!-\!2k\big\}
\sbnm{n}{k}
\frac{\binm{x+n}{k}
\binm{n+bn}{k}
\binm{n+dn}{k}}
{\binm{k-x}{k}
\binm{k+bn-x}{k}
\binm{k+dn-x}{k}}
=(-1)^nx
\frac{\binm{x+n}{n}
\binm{2n+bn+dn-x}{n}}
{\binm{n+bn-x}{n}
\binm{n+dn-x}{n}}\]
which leads us, under the derivative operator 
$\mc{D}_0$, to the following result.
\begin{thm}\label{thm:0/2} With $b,\:d\in\mb{N}_0$, 
there holds the following harmonic number identity:
\[\sum_{k=0}^n
\sbnm{n}{k}^2
\frac{\binm{n+bn}{k}
\binm{n+dn}{k}}
{\binm{k+bn}{k}
\binm{k+dn}{k}}
\Big\{1+(n\!-\!2k)
\big(2H_k\!+\!H_{bn+k}
\!+\!H_{dn+k}\big)\Big\}
=(-1)^n
\frac{\binm{2n+bn+dn}{n}}
{\binm{n+bn}{n}
\binm{n+dn}{n}}.\]
\end{thm}


\section{The Whipple Transformation}
In this section, the Whipple transformation~\citu{bailey}{\S4.3} 
\xalignz{
_{7}F_{6}\ffnk{rccccc}{\:1}
{a,\;1+a/2,\+b,\+c,\+d,\+e,\+-n}
{a/2,\+1+a-b,\+1+a-c,\+1+a-d,
\+1+a-e,\+1+a+n}&\\[2.5mm]
=\frac{(1\!+\!a)_n(1\!+\!a\!-\!b\!-\!d)_n}
      {(1\!+\!a\!-\!b)_n(1\!+\!a\!-\!d)_n}
{_{4}F_{3}\ffnk{cccc}{1}
{-n,\+b,\+d,\+1+a-c-e}
{\+1+a-c,\+1+a-e,\+b+d-a-n}}
&}
will be used to derive harmonic number identities.


\subsection{}
Performing parameter replacement   
\[\left.\xalignedx{
a&\to-x-n\\
b&\to1+bn\\
c&\to1+cn\\
d&\to1+dn\\
e&\to1+en}
\right\}\qquad
(b,\:c,\:d,\:e\in\mb{N}_0)\]
we can restate the Whipple transformation as
\xalignz{
&\sum_{k=0}^n
\big\{x\!+\!n\!-\!2k\big\}
\sbnm{n}{k}
\frac{\binm{x+n}{k}
\binm{k+bn}{k}
\binm{k+cn}{k}
\binm{k+dn}{k}
\binm{k+en}{k}}
{\binm{k-x}{k}
\binm{n+bn+x}{k}
\binm{n+cn+x}{k}
\binm{n+dn+x}{k}
\binm{n+en+x}{k}}\\
&\:=\:x
\frac{\binm{x+n}{n}
\binm{1+x+bn+dn+n}{n}}
{\binm{x+bn+n}{n}
\binm{x+dn+n}{n}}
\sum_{\ell=0}^n
\sbnm{n}{\ell}
\frac{\binm{\ell+bn}{\ell}
\binm{\ell+dn}{\ell}
\binm{1+x+cn+en+n}{\ell}}
{\binm{x+cn+n}{\ell}
\binm{x+en+n}{\ell}
\binm{1+x+bn+dn+\ell}{\ell}}}
which leads us, under the derivative operator 
$\mc{D}_0$, to the following result.
\begin{thm}\label{thm:4/0} 
For four nonnegative integers 
$\{b,\:c,\:d,\:e\}$, there holds:
\xalignz{
&\sum_{k=0}^n
\sbnm{n}{k}^2
\frac{\binm{k+bn}{k}
\binm{k+cn}{k}
\binm{k+dn}{k}
\binm{k+en}{k}}
{\binm{n+bn}{k}
\binm{n+cn}{k}
\binm{n+dn}{k}
\binm{n+en}{k}}\\
&\times\Big\{1\!+\!(n\!-\!2k)
\big(2H_k\!-\!H_{bn+k}
\!-\!H_{cn+k}\!-\!H_{dn+k}
\!-\!H_{en+k}\big)\Big\}\\
&\:=\:
\frac{\binm{1+bn+dn+n}{n}}
{\binm{n+bn}{n}\binm{n+dn}{n}}
\sum_{\ell=0}^n
\sbnm{n}{\ell}
\frac{\binm{\ell+bn}{\ell}
\binm{\ell+dn}{\ell}
\binm{1+cn+en+n}{\ell}}
{\binm{n+cn}{\ell}
\binm{n+en}{\ell}
\binm{1+bn+dn+\ell}{\ell}}.}
\end{thm}


\subsection{}
Performing parameter replacement   
\[\left.\xalignedx{
a&\to-x-n\\
b&\to1+bn\\
c&\to1+cn\\
d&\to1+dn\\
e&\to-n-en}
\right\}\qquad
(b,\:c,\:d,\:e\in\mb{N}_0)\]
we can restate the Whipple transformation as
\xalignz{
&\sum_{k=0}^n
\big\{x\!+\!n\!-\!2k\big\}
\sbnm{n}{k}
\frac{\binm{x+n}{k}
\binm{k+bn}{k}
\binm{k+cn}{k}
\binm{k+dn}{k}
\binm{n+en}{k}}
{\binm{k-x}{k}
\binm{n+bn+x}{k}
\binm{n+cn+x}{k}
\binm{n+dn+x}{k}
\binm{k+en-x}{k}}\\
&\:=\:x
\frac{\binm{x+n}{n}
\binm{1+x+bn+dn+n}{n}}
{\binm{x+bn+n}{n}
\binm{x+dn+n}{n}}
\sum_{\ell=0}^n(-1)^{\ell}
\sbnm{n}{\ell}
\frac{\binm{\ell+bn}{\ell}
\binm{\ell+dn}{\ell}
\binm{x+cn-en}{\ell}}
{\binm{n+cn+x}{\ell}
\binm{\ell+en-x}{\ell}
\binm{1+x+bn+dn+\ell}{\ell}}}
which leads us, under the derivative operator 
$\mc{D}_0$, to the following result.
\begin{thm}\label{thm:3/1}
For four nonnegative integers 
$\{b,\:c,\:d,\:e\}$, there holds:
\xalignz{
&\sum_{k=0}^n
\sbnm{n}{k}^2
\frac{\binm{k+bn}{k}
\binm{k+cn}{k}
\binm{k+dn}{k}
\binm{n+en}{k}}
{\binm{n+bn}{k}
\binm{n+cn}{k}
\binm{n+dn}{k}
\binm{k+en}{k}}\\
&\times\Big\{1\!+\!(n\!-\!2k)
\big(2H_k\!-\!H_{bn+k}
\!-\!H_{cn+k}\!-\!H_{dn+k}
\!+\!H_{en+k}\big)\Big\}\\
&\:=\:
\frac{\binm{1+bn+dn+n}{n}}
{\binm{n+bn}{n}\binm{n+dn}{n}}
\sum_{\ell=0}^n(-1)^{\ell}
\sbnm{n}{\ell}
\frac{\binm{\ell+bn}{\ell}
\binm{\ell+dn}{\ell}
\binm{cn-en}{\ell}}
{\binm{n+cn}{\ell}
\binm{\ell+en}{\ell}
\binm{1+bn+dn+\ell}{\ell}}.}
\end{thm}


\subsection{}
Performing parameter replacement   
\[\left.\xalignedx{
a&\to-x-n\\
b&\to1+bn\\
c&\to-n-cn\\
d&\to1+dn\\
e&\to-n-en}
\right\}\qquad
(b,\:c,\:d,\:e\in\mb{N}_0)\]
we can restate the Whipple transformation as
\xalignz{
&\sum_{k=0}^n
\big\{x\!+\!n\!-\!2k\big\}
\sbnm{n}{k}
\frac{\binm{x+n}{k}
\binm{k+bn}{k}
\binm{n+cn}{k}
\binm{k+dn}{k}
\binm{n+en}{k}}
{\binm{k-x}{k}
\binm{n+bn+x}{k}
\binm{k+cn-x}{k}
\binm{n+dn+x}{k}
\binm{k+en-x}{k}}\\
&\:=\:x
\frac{\binm{x+n}{n}
\binm{1+x+bn+dn+n}{n}}
{\binm{x+bn+n}{n}
\binm{x+dn+n}{n}}
\sum_{\ell=0}^n(-1)^{\ell}
\sbnm{n}{\ell}
\frac{\binm{\ell+bn}{\ell}
\binm{\ell+dn}{\ell}
\binm{\ell+n+cn+en-x}{\ell}}
{\binm{\ell+cn-x}{\ell}
\binm{\ell+en-x}{\ell}
\binm{1+x+bn+dn+\ell}{\ell}}}
which leads us, under the derivative operator 
$\mc{D}_0$, to the following result.
\begin{thm}\label{thm:2/2}
For four nonnegative integers 
$\{b,\:c,\:d,\:e\}$, there holds:
\xalignz{
&\sum_{k=0}^n
\sbnm{n}{k}^2
\frac{\binm{k+bn}{k}
\binm{n+cn}{k}
\binm{k+dn}{k}
\binm{n+en}{k}}
{\binm{n+bn}{k}
\binm{k+cn}{k}
\binm{n+dn}{k}
\binm{k+en}{k}}\\
&\times\Big\{1\!+\!(n\!-\!2k)
\big(2H_k\!-\!H_{bn+k}
\!+\!H_{cn+k}\!-\!H_{dn+k}
\!+\!H_{en+k}\big)\Big\}\\
&\:=\:
\frac{\binm{1+bn+dn+n}{n}}
{\binm{n+bn}{n}\binm{n+dn}{n}}
\sum_{\ell=0}^n(-1)^{\ell}
\sbnm{n}{\ell}
\frac{\binm{\ell+bn}{\ell}
\binm{\ell+dn}{\ell}
\binm{n+cn+en+\ell}{\ell}}
{\binm{\ell+cn}{\ell}
\binm{\ell+en}{\ell}
\binm{1+bn+dn+\ell}{\ell}}.}
\end{thm}


\subsection{}
Performing parameter replacement   
\[\left.\xalignedx{
a&\to-x-n\\
b&\to1+bn\\
c&\to-n-cn\\
d&\to-n-dn\\
e&\to-n-en}
\right\}\qquad
(b,\:c,\:d,\:e\in\mb{N}_0)\]
we can restate the Whipple transformation as
\xalignz{
&\sum_{k=0}^n
\big\{x\!+\!n\!-\!2k\big\}
\sbnm{n}{k}
\frac{\binm{x+n}{k}
\binm{k+bn}{k}
\binm{n+cn}{k}
\binm{n+dn}{k}
\binm{n+en}{k}}
{\binm{k-x}{k}
\binm{n+bn+x}{k}
\binm{k+cn-x}{k}
\binm{k+dn-x}{k}
\binm{k+en-x}{k}}\\
&\:=\:(-1)^nx
\frac{\binm{x+n}{n}
\binm{bn-dn+x}{n}}
{\binm{x+bn+n}{n}
\binm{n+dn-x}{n}}
\sum_{\ell=0}^n
\sbnm{n}{\ell}
\frac{\binm{\ell+bn}{\ell}
\binm{n+dn}{\ell}
\binm{\ell+n+cn+en-x}{\ell}}
{\binm{\ell+cn-x}{\ell}
\binm{\ell+en-x}{\ell}
\binm{\ell-n+bn-dn+x}{\ell}}}
which leads us, under the derivative operator 
$\mc{D}_0$, to the following result.
\begin{thm}\label{thm:1/3}
For four nonnegative integers 
$\{b,\:c,\:d,\:e\}$, there holds:
\xalignz{
&\sum_{k=0}^n
\sbnm{n}{k}^2
\frac{\binm{k+bn}{k}
\binm{n+cn}{k}
\binm{n+dn}{k}
\binm{n+en}{k}}
{\binm{n+bn}{k}
\binm{k+cn}{k}
\binm{k+dn}{k}
\binm{k+en}{k}}\\
&\times\Big\{1\!+\!(n\!-\!2k)
\big(2H_k\!-\!H_{bn+k}
\!+\!H_{cn+k}\!+\!H_{dn+k}
\!+\!H_{en+k}\big)\Big\}\\
&\:=\:(-1)^n
\frac{\binm{bn-dn}{n}}
{\binm{n+bn}{n}\binm{n+dn}{n}}
\sum_{\ell=0}^n
\sbnm{n}{\ell}
\frac{\binm{\ell+bn}{\ell}
\binm{n+dn}{\ell}
\binm{\ell+cn+en+n}{\ell}}
{\binm{\ell+cn}{\ell}
\binm{\ell+en}{\ell}
\binm{\ell+bn-dn-n}{\ell}}.}
\end{thm}


\subsection{}
Performing parameter replacement   
\[\left.\xalignedx{
a&\to-x-n\\
b&\to-n-bn\\
c&\to-n-cn\\
d&\to-n-dn\\
e&\to-n-en}
\right\}\qquad
(b,\:c,\:d,\:e\in\mb{N}_0)\]
we can restate the Whipple transformation as
\xalignz{
&\sum_{k=0}^n
\big\{x\!+\!n\!-\!2k\big\}
\sbnm{n}{k}
\frac{\binm{x+n}{k}
\binm{n+bn}{k}
\binm{n+cn}{k}
\binm{n+dn}{k}
\binm{n+en}{k}}
{\binm{k-x}{k}
\binm{k+bn-x}{k}
\binm{k+cn-x}{k}
\binm{k+dn-x}{k}
\binm{k+en-x}{k}}\\
&\:=\:(-1)^nx
\frac{\binm{x+n}{n}
\binm{2n+bn+dn-x}{n}}
{\binm{n+bn-x}{n}
\binm{n+dn-x}{n}}
\sum_{\ell=0}^n
\sbnm{n}{\ell}
\frac{\binm{n+bn}{\ell}
\binm{n+dn}{\ell}
\binm{\ell+n+cn+en-x}{\ell}}
{\binm{\ell+cn-x}{\ell}
\binm{\ell+en-x}{\ell}
\binm{2n+bn+dn-x}{\ell}}}
which leads us, under the derivative operator 
$\mc{D}_0$, to the following result.
\begin{thm}\label{thm:0/4}
For four nonnegative integers 
$\{b,\:c,\:d,\:e\}$, there holds:
\xalignz{
&\sum_{k=0}^n
\sbnm{n}{k}^2
\frac{\binm{n+bn}{k}
\binm{n+cn}{k}
\binm{n+dn}{k}
\binm{n+en}{k}}
{\binm{k+bn}{k}
\binm{k+cn}{k}
\binm{k+dn}{k}
\binm{k+en}{k}}\\
&\times\Big\{1\!+\!(n\!-\!2k)
\big(2H_k\!+\!H_{bn+k}
\!+\!H_{cn+k}\!+\!H_{dn+k}
\!+\!H_{en+k}\big)\Big\}\\
&\:=\:(-1)^n
\frac{\binm{2n+bn+dn}{n}}
{\binm{n+bn}{n}\binm{n+dn}{n}}
\sum_{\ell=0}^n
\sbnm{n}{\ell}
\frac{\binm{n+bn}{\ell}
\binm{n+dn}{\ell}
\binm{n+cn+en+\ell}{\ell}}
{\binm{\ell+cn}{\ell}
\binm{\ell+en}{\ell}
\binm{2n+bn+dn}{\ell}}.}
\end{thm}


\section{Harmonic Number Identities and Transformations}
\setcounter{equation}{0}

In order to facilitate computation of harmonic number sums,
we present a useful limiting relation concerning harmonic numbers. 
Suppose that $\lam,\nu,n,k\in\mb{N}_0$ (the set of nonnegative integers)
with $k\le n$ and $\{P_k(y),\:Q_k(y)\}$ are two families of monic 
polynomials with $P_k(y)$ and $Q_k(y)$ being of degree $k$ in $y$, 
then there holds 
\eqn{\lim_{y\to\infty}\Big\{
\frac{P_{\lam k+\nu}(y)}
     {Q_{\lam k+\nu}(y)}
H_{ny+k}\label{diff}
-\frac{P_{\lam(n-k)+\nu}(y)}
      {Q_{\lam(n-k)+\nu}(y)}
H_{ny+n-k}\bigg\}=0.} 

In fact, it is not hard to see that 
$\tfrac{P_{\lam k+\nu}(y)}{Q_{\lam k+\nu}(y)}$
tends to one and $H_{ny+k}\approx\ln(ny+k)$
as $y\to\infty$. Now rewrite the function
in question into two terms 
\bnm
\frac{P_{\lam k+\nu}(y)}
     {Q_{\lam k+\nu}(y)}
H_{ny+k}
-\frac{P_{\lam(n-k)+\nu}(y)}
      {Q_{\lam(n-k)+\nu}(y)}
H_{ny+n-k}
&=&\Big\{ H_{ny+k}-H_{ny+n-k}\Big\}
\frac{P_{\lam(n-k)+\nu}(y)}
     {Q_{\lam(n-k)+\nu}(y)}\\
&+&H_{ny+k}\bigg\{ 
\frac{P_{\lam k+\nu}(y)}
     {Q_{\lam k+\nu}(y)}
-\frac{P_{\lam(n-k)+\nu}(y)}
      {Q_{\lam(n-k)+\nu}(y)}\bigg\}.
\enm
When $y\to\infty$, the right hand side on the penultimate line
tends to zero because the fraction is bounded and the difference  
in braces behaves like $\ln\frac{ny+k}{ny+n-k}\to0$; 
the last line tends to zero too since the fractional difference 
is a fraction with numerator degree less than denominator
degree in view of the fact that both $P(y)$ and $Q(y)$ are 
polynomials with the leading coefficients equal to one. 

\begin{thm}\label{thm:limit}
Let $\{P_k(y),\:Q_k(y)\}$ be two families of monic 
polynomials with $P_k(y)$ and $Q_k(y)$ being of degree 
$k$ in $y$. If $f_n(k)$ is a function independent of $y$ 
which satisfies the reflection property $f_n(k)=-f_n(n-k)$,
then there holds the following limiting relation: 
\eqn{\lim_{y\to\infty}
\sum_{k=0}^nf_n(k)
\frac{P_{\lam k+\nu}(y)}
     {Q_{\lam k+\nu}(y)}
H_{ny+k}=0.\label{limit}} 
\end{thm}
\begin{proof}
By means of the summation index involution 
$k\to n-k$, we can reformulate the finite 
sum stated in the theorem as 
\[\sum_{k=0}^nf_n(k)
\frac{P_{\lam k+\nu}(y)}
     {Q_{\lam k+\nu}(y)}
H_{ny+k}
\:=\:\frac{1}{2}
\sum_{k=0}^nf_n(k)
\bigg\{\frac{P_{\lam k+\nu}(y)}
     {Q_{\lam k+\nu}(y)}
H_{ny+k}
-\frac{P_{\lam(n-k)+\nu}(y)}
      {Q_{\lam(n-k)+\nu}(y)}
H_{ny+n-k}\bigg\}.\] 
In view of \eqref{diff}, the differences in the braces 
on the right hand side tends to zero as $y\to\infty$. 
We therefore obtain the limiting relation about 
harmonic number sums stated in the theorem.
\end{proof}

There is a large class of functions satisfying 
the reflection property in the theorem, for example 
\eqn{f_n(k)\:=\:\sbnm{n}{k}^\mu
\frac{\binm{n+k}{k}^\nu}{\binm{2n}{k}^\nu}
(n-2k),\qquad(\mu,\:\nu\in\mb{N}_0)\label{reflex}}
which come out frequently for the limiting process
in the construction of Table-I and Table-II.


Now we take Entry-4 from Table-II to exemplify  
how to derive harmonic number identities 
from the theorems established in this paper.
 
Specifying with $b=d=1$ and $e=0$,
we can state the transformation in Theorem~\ref{thm:3/1}
as
\[\sum_{k=0}^n\sbnm{n}{k}^3
\frac{\binm{n+k}{k}^2}{\binm{2n}{k}^2}
\frac{\binm{k+cn}{k}}{\binm{n+cn}{k}}
\Big\{{\sst1\!+\!(n\!-\!2k)
\big(3H_k-2H_{n+k}-H_{cn+k}\big)}\Big\}
\:=\:
\frac{\binm{1+3n}{n}}{\binm{2n}{n}^2}
\sum_{\ell=0}^n(-1)^{\ell}
\frac{\binm{n}{\ell}\binm{n+\ell}{\ell}^2\binm{cn}{\ell}}
     {\binm{1+2n+\ell}{\ell}\binm{n+cn}{\ell}}.\]
It is easy to see that the coefficient corresponding 
to $H_{cn+k}$ is given by \eqref{reflex} with $\mu=3$ 
and $\nu=2$. In view of Theorem~\ref{thm:limit}, 
the limit $c\to\infty$ of the last equation reads as
\eqn{\sum_{k=0}^n\sbnm{n}{k}^3
\frac{\binm{n+k}{k}^2}{\binm{2n}{k}^2}
\Big\{1\!+\!(n\!-\!2k)\big(3H_k-2H_{n+k}\big)\Big\}
\:=\:\frac{\binm{1+3n}{n}}{\binm{2n}{n}^2}
\sum_{\ell=0}^n(-1)^{\ell}
\frac{\binm{n}{\ell}\binm{n+\ell}{\ell}^2}
     {\binm{1+2n+\ell}{\ell}}}
which is exactly the fourth identity displayed in Table-II.

We remark that the right hand side of this last identity 
can further be evaluated by Dixon's formula and 
we therefore get the following closed formula:
\eqn{\sum_{k=0}^n\sbnm{n}{k}^3
\frac{\binm{n+k}{k}^2}{\binm{2n}{k}^2}
\Big\{1\!+\!(n\!-\!2k)\big(3H_k-2H_{n+k}\big)\Big\}
\:=\:\begin{cases}
0,&n-\text{odd}\\
(-1)^{m}\binm{3m}{m,m,m}
\bigg/\binm{4m}{2m}^2,
&n=2m.\end{cases}}

 
Specifying the free parameters in 
Theorems~\ref{thm:ps-lam}-\ref{thm:0/4}, 
we can similarly establish, by means of Theorem~\ref{thm:limit}, 
$26$ closed summation formulas and $21$ transformations 
on harmonic numbers, which are displayed respectively 
in Table-I and Table-II at the end of this paper.


As a partial answer to the question posed at the end of the paper
by Paule and Schneider~\cito{paule}, the examples~8,\:9,\:16,\:17
numbered with in Table-I and 16\,\:17 in Table-II confirm that the sum
\[\Xi_\lam(n)\::=\:\sum_{k=0}^n
\sbnm{n}{k}^\lam 
\Big\{1+\lam(n-2k)\:H_k\Big\}
\qquad(\lam,\:n\in\mb{N})\]
are representable in terms of terminating hypergeometric series
for $1\le\lam\le6$. In addition, the hypergeometric method
presented in this paper shows that these binomial-harmonic number 
sums trace back to the same origin - the very-well poised
terminating hypergeometric series. In fact, if we define
\[\Omega_\lam(n,x)\::=\:{_{1+\lam}F_\lam}
\ffnk{crc}{(-1)^\lam}
  {-x-n,&1-\frac{x+n}{2},&\ang{-n}_{\lam-1}}
  {&-\frac{x+n}{2},&\ang{1-x}_{\lam-1}\rule{0mm}{5mm}}\]
where $\ang{w}_{\lam}$ stands for $\lam$ copies of $w$.
Then it is not difficult to check that
\[\Xi_{\lam}(n)\:=\:\mc{D}_0\Big\{(x+n)\Omega_\lam(n,x)\Big\}.\]
However, the problem posed by Paule and Schneider~\cito{paule}
remains open for $\lam>6$, i.e., whether $\Xi_\lam(n)$ can
be expressed as a single terminating hypergeometric series.

\thank{145mm}{In a recent preprint 
``Hyperg\éom\étrie et fonction z\^eta de Riemann" 
by Christian Krattenthaler and Tanguy Rivoal,
a multisum expression for $\Xi_\lam(n)$ has been 
derived, but as pointed out by Krattenthaler
to the authors, that it is \emph{(most likely)} not 
possible to express these sums as single hypergeometric sums. 
They make also the same observation, namely that the
identities proved in the paper by Paule and Schneider 
come from applying differentiation to known hypergeometric 
summation or transformation theorems. In this sense, their 
work has some common background with ours, 
but they have different aims. 
The authors thank to Krattenthaler for the information.}

\newpage

\centro{\textbf{Table - I}: 
The harmonic number identities of type 
$\dst\sum_{k=0}^nA(n,k)=C(n)$:
\nopagebreak[4]
\begin{tabular}{||c|l|l|c||}
\hline\hline
\makebox[1mm]{No}
&\qquad$\xquad A(n,k)$
&\:$\xquad C(n)$
&Note\rule[-1.5mm]{0mm}{5mm}\\
\hline\hline1\rule[-2.5mm]{0mm}{7mm}
&$\binm{n}{k}^2\binm{2n+k}{k}
\big\{\sst H_{2n+k}-H_{k}\big\}$
&$2\binm{2n}{n}^2
\big\{\sst H_{2n}-H_{n}\big\}$
&$\double{\text{Thm\,\ref{thm:ps-lam}:}
	\:\lam=2}{\mu=\nu=0}$\\
\hline2\rule[-2.5mm]{0mm}{7mm}
&$\binm{n}{k}^2\binm{2n+k}{k}
\:\big\{\sst H_{k}-H_{n-k}\big\}$
&$\binm{2n}{n}^2
\big\{\sst H_{2n}-H_{n}\big\}$
&$\double{\text{Thm\,\ref{thm:ps-mu}:}
	\:\lam=2}{\mu=\nu=0}$\\
\hline3\rule[-2.5mm]{0mm}{7mm}
&$\binm{n}{k}\binm{2n}{k}\binm{3n+k}{k}
\big\{\sst H_{3n+k}-H_{k}\big\}$
&$\binm{3n}{n}^2
\big\{\sst2H_{3n}-H_{n}-H_{2n}\big\}$
&$\double{\text{Thm\,\ref{thm:ps-lam}:}
	\:\lam=3}{\mu=1\text{ and }\nu=0}$\\
\hline4\rule[-2.5mm]{0mm}{7mm}
&$\binm{n}{k}
\binm{2n}{k}\binm{3n+k}{k}
\big\{\sst H_{2n-k}-H_{k}\big\}$
&$\binm{3n}{n}^2
\big\{\sst2H_{2n}-H_{n}-H_{3n}\big\}$
&$\double{\text{Thm\,\ref{thm:ps-mu}:}
	\:\lam=3}{\mu=1\text{ and }\nu=0}$\\
\hline5\rule[-2.5mm]{0mm}{7mm}
&$\binm{n}{k}^2
\binm{3n+k}{2n}
\big\{\sst H_{3n+k}-H_{k}\big\}$
&$\binm{3n}{n}
\big\{\sst H_{2n}+H_{3n}-2H_{n}\big\}$
&$\double{\text{Thm\,\ref{thm:ps-lam}:}
	\:\lam=3}{\mu=0\text{ and }\nu=1}$\\
\hline6\rule[-2.5mm]{0mm}{7mm}
&$\binm{n}{k}^2
\binm{3n+k}{2n}
\big\{\sst H_{k}-H_{n-k}\big\}$
&$\binm{3n}{n}
\big\{\sst H_{3n}-H_{2n}\big\}$
&$\double{\text{Thm\,\ref{thm:ps-mu}:}
	\:\lam=3}{\mu=0\text{ and }\nu=1}$\\
\hline7\rule[-2.5mm]{0mm}{7mm}
&$\binm{n}{k}^2
\binm{3n+k}{2n}
\big\{\sst H_{n+k}-H_{k}\big\}$
&$\binm{3n}{n}
\big\{\sst3H_{2n}-2H_{n}-H_{3n}\big\}$
&$\double{\text{Thm\,\ref{thm:ps-nu}:}
	\:\lam=3}{\mu=0\text{ and }\nu=1}$\\
\hline8\rule[-2.5mm]{0mm}{7mm}
&$\binm{n}{k}
\big\{\sst1+(n-2k)H_k\big\}$
&$1$
&$\double{\text{Thm\,\ref{thm:2/0}:}\:b=0}
         {d\to\infty;\:\text{cf.\:\citu{paule}{Eq\:1}}}$\\
\hline9\rule[-2.5mm]{0mm}{7mm}
&$\binm{n}{k}^2
\big\{\sst1+2(n-2k)H_k\big\}$
&$0$
&$\double{\text{Thm\,\ref{thm:2/0}:}\:b,\:d\to\infty}
	   {\text{cf.\:\citu{paule}{Eq\:2}}}$\\
\hline10\rule[-2.5mm]{0mm}{7mm}
&$\binm{n+k}{k}\binm{2n-k}{n}
\big\{\sst1+(n-2k)(H_k-H_{n+k})\big\}$
&$\binm{1+2n}{n}$
&$\sst{\text{Thm\,\ref{thm:2/0}:}
	\:b=0\:\&\:d=1}{}$\\
\hline11\rule[-2.5mm]{0mm}{7mm}
&$\binm{n+k}{k}^2\binm{2n-k}{n}^2
\big\{\sst1+2(n-2k)(H_k-H_{n+k})\big\}$
&$\binm{1+3n}{n}$
&$\sst{\text{Thm\,\ref{thm:2/0}:}
	\:b=d=1}{}$\\
\hline12\rule[-2.5mm]{0mm}{7mm}
&$\binm{2n}{k}\binm{2n}{n+k}
\big\{\sst1+(n-2k)(H_k+H_{n+k})\big\}$
&$\binm{2n-1}{n}$
&$\sst{\text{Thm\,\ref{thm:1/1}:}
	\:b=0\:\&\:d=1}{}$\\
\hline13\rule[-2.5mm]{0mm}{7mm}
&$\binm{n}{k}
\binm{2n}{k}\binm{2n}{n+k}
\big\{\sst1+(n-2k)(2H_k+H_{n+k})\big\}$
&$(-1)^n$
&$\sst{\text{Thm\,\ref{thm:1/1}:}
	\:b\to\infty\:\&\:d=1}{}$\\
\hline14\rule[-2.5mm]{0mm}{7mm}
&$\binm{n}{k}
\binm{n+k}{n}\binm{2n-k}{n}
\big\{\sst1+(n-2k)(2H_k-H_{n+k})\big\}$
&$1$
&$\sst{\text{Thm\,\ref{thm:1/1}:}
	\:b=1\:\&\:d\to\infty}{}$\\
\hline15\rule[-2.5mm]{0mm}{7mm}
&$\binm{n}{k}^2
\binm{n+k}{n}\binm{2n-k}{n}
\big\{\sst1+(n-2k)(3H_k-H_{n+k})\big\}$
&$(-1)^n$
&$\sst{\text{Thm\,\ref{thm:1/1}:}
	\:b=1\:\&\:d=0}{}$\\
\hline16\rule[-2.5mm]{0mm}{7mm}
&$\binm{n}{k}^3
\big\{\sst1+3(n-2k)H_k\big\}$
&$(-1)^n$
&$\double{\text{Thm\,\ref{thm:0/2}:}b=0}
	   {d\to\infty;\:\text{cf.\:\citu{paule}{Eq\:3}}}$\\
\hline17\rule[-2.5mm]{0mm}{7mm}
&$\binm{n}{k}^4
\big\{\sst1+4(n-2k)H_k\big\}$
&$(-1)^{n}\binm{2n}{n}$
&$\double{\text{Thm\,\ref{thm:0/2}:\:}b=d=0}
	   {\text{cf.\:\citu{paule}{Eq\:4}}}$\\
\hline18\rule[-2.5mm]{0mm}{7mm}
&$\binm{n}{k}^2\binm{2n}{k}\binm{2n}{n+k}
\big\{\sst1+(n-2k)(3H_k+H_{n+k})\big\}$
&$(-1)^{n}
\binm{3n}{n}$
&$\sst{\text{Thm\,\ref{thm:0/2}:}
	\:b=0\:\&\:d=1}{}$\\
\hline19\rule[-2.5mm]{0mm}{7mm}
&$\binm{2n}{k}^2\binm{2n}{n+k}^2
\big\{\sst1+2(n-2k)(H_k+H_{n+k})\big\}$
&$(-1)^{n}
\binm{4n}{n}$
&$\sst{\text{Thm\,\ref{thm:0/2}:}
	\:b=d=1}{}$\\
\hline20\rule[-2.5mm]{0mm}{7mm}
&$\binm{n}{k}^{-1}
\big\{\sst1-(n-2k)H_k\big\}$
&$(1+n)H_{n+1}$
&$\double{\text{Thm\,\ref{thm:4/0}:}
	\:e\to\infty}{b=c=d=0}$\\
\hline21\rule[-2.5mm]{0mm}{7mm}
&$\binm{n}{k}^{-2}
\big\{\sst1-2(n-2k)H_k\big\}$
&$2\frac{(1+n)^2}{2+n}H_{n+1}$
&$\double{\text{Thm\,\ref{thm:4/0}:}
	\:}{b=c=d=e=0}$\\
\hline22\rule[-3.5mm]{0mm}{9mm}
&$\frac{1-(n-2k)(H_k+H_{n+k})}
     {\binm{2n}{k}\binm{2n}{n+k}}$
&$\frac{1+2n}{2+2n}
+(\sst n+\frac{1}{2})H_{1+2n}$
&$\double{\text{Thm\,\ref{thm:4/0}:}
	\:e=1}{b=c=d=0}$\\
\hline23\rule[-3.5mm]{0mm}{9mm}
&$\frac{\binm{n+k}{k}}{\binm{2n}{k}}
\big\{\sst1-(n-2k)H_{n+k}\big\}$
&${(1+2n)}\big\{\sst H_{1+2n}-H_{n}\big\}$
&$\double{\text{Thm\,\ref{thm:4/0}:}
	\:b=d=0}{c=1\:\&\:e\to\infty}$\\
\hline24\rule[-3.5mm]{0mm}{9mm}
&$\frac{\binm{n+k}{k}^2}{\binm{2n}{k}^2}
\big\{\sst1-2(n-2k)H_{n+k}\big\}$
&$2\frac{(1+2n)^2}{2+3n}
\big\{\sst H_{1+2n}-H_{n}\big\}$
&$\double{\text{Thm\,\ref{thm:4/0}:}
	\:b=d=0}{c=e=1}$\\
\hline25\rule[-3.5mm]{0mm}{9mm}
&$\frac{\binm{2n}{k}}{\binm{n}{k}\binm{n+k}{k}}
\big\{\sst1-(n-2k)(H_k-H_{n+k})\big\}$
&$\frac{n(n+1)}{n-1}
\big\{\sst H_{n+1}+H_{n-1}-H_{2n}\big\}$
&$\double{\text{Thm\,\ref{thm:3/1}:}
	\:n>1}{b=c=d=0\:\&\:e=1}$\\
\hline26\rule[-3.5mm]{0mm}{9mm}
&$\frac{\binm{2n}{k}^2}{\binm{n+k}{k}^2}
\big\{\sst1+2(n-2k)H_{n+k}\big\}$
&$\frac{2n}{3}
\big\{\sst H_{2n}-H_{n-1}\big\}$
&$\double{\text{Thm\,\ref{thm:2/2}:}
	\:n>0}{b=d=0\:\&\:c=e=1}$\\
\hline\hline
\end{tabular}}


\centro{\textbf{Table - II}: 
The harmonic number transformations of type 
$\dst\sum_{k=0}^nA(n,k)=\dst\sum_{\ell=0}^nB(n,\ell)$:
\nopagebreak[4]
\begin{tabular}{||c|l|l|c||}
\hline\hline
\makebox[1mm]{No}
&\qquad$\xquad A(n,k)$
&\quad$\xquad B(n,\ell)$
&Note\rule[-1.5mm]{0mm}{5mm}\\
\hline1\rule[-3.5mm]{0mm}{9mm}
&$\binm{n}{k}
\frac{\binm{n+k}{k}^3}{\binm{2n}{k}^3}
\big\{\sst1+(n-2k)(H_k-3H_{n+k})\big\}$
&$\frac{\binm{1+3n}{n}}{\binm{2n}{n}^2}
\times\frac{1+2n}{1+2n-\ell}
\frac{\binm{n+\ell}{\ell}^2}{\binm{1+2n+\ell}{\ell}}$
&$\double{\text{Thm\,\ref{thm:4/0}:}
	\:e=0}{b=c=d=1}$\\
\hline2\rule[-3.5mm]{0mm}{9mm}
&$\binm{n}{k}^2
\frac{\binm{n+k}{k}^4}{\binm{2n}{k}^4}
\big\{\sst1+2(n-2k)(H_k-2H_{n+k})\big\}$
&$\frac{\binm{1+3n}{n}}{\binm{2n}{n}^2}
\times\binm{n}{\ell}
\frac{\binm{n+\ell}{\ell}^2\binm{1+3n}{\ell}}
     {\binm{2n}{\ell}^2\binm{1+2n+\ell}{\ell}}$
&$\double{\text{Thm\,\ref{thm:4/0}:}}
	   {b=c=d=e=1}$\\
\hline3\rule[-3.5mm]{0mm}{9mm}
&$\binm{n}{k}
\frac{\binm{n+k}{k}^2}{\binm{2n}{k}^2}
\big\{\sst1+(n-2k)(H_k-2H_{n+k})\big\}$
&$\frac{\binm{1+3n}{n}}{\binm{2n}{n}^2}
\times
\frac{\binm{n+\ell}{\ell}^2}{\binm{1+2n+\ell}{\ell}}$
&$\double{\text{Thm\,\ref{thm:3/1}:}
	\:b=d=1}{c=0\:\&\:e\to\infty}$\\
\hline4\rule[-3.5mm]{0mm}{9mm}
&$\binm{n}{k}^3
\frac{\binm{n+k}{k}^2}{\binm{2n}{k}^2}
\big\{\sst1+(n-2k)(3H_k-2H_{n+k})\big\}$
&$\frac{\binm{1+3n}{n}}{\binm{2n}{n}^2}
\times(-1)^{\ell}
\frac{\binm{n}{\ell}\binm{n+\ell}{\ell}^2}
     {\binm{1+2n+\ell}{\ell}}$
&$\double{\text{Thm\,\ref{thm:3/1}:}
	\:b=d=1}{c\to\infty\:\&\:e=0}$\\
\hline5\rule[-3.5mm]{0mm}{9mm}
&$\binm{n}{k}^2
\frac{\binm{n+k}{k}^3}{\binm{2n}{k}^3}
\big\{\sst1+(n-2k)(2H_k-3H_{n+k})\big\}$
&$\frac{\binm{1+3n}{n}}{\binm{2n}{n}^2}
\times
\frac{\binm{n}{\ell}\binm{n+\ell}{\ell}^2}
     {\binm{2n}{\ell}\binm{1+2n+\ell}{\ell}}$
&$\double{\text{Thm\,\ref{thm:3/1}:}
	\:e\to\infty}{b=c=d=1}$\\
\hline6\rule[-3.5mm]{0mm}{9mm}
&$\binm{n}{k}^3\frac{\binm{n+k}{k}^3}{\binm{2n}{k}^3}
\big\{\sst1+3(n-2k)(H_k-H_{n+k})\big\}$
&$\frac{\binm{1+3n}{n}}{\binm{2n}{n}^2}
\times(-1)^{\ell}
\frac{\binm{n}{\ell}^2\binm{n+\ell}{\ell}^2}
{\binm{2n}{\ell}\binm{1+2n+\ell}{\ell}}$
&$\double{\text{Thm\:\ref{thm:3/1}:}
	\:e=0}{b=c=d=1}$\\
\hline7\rule[-3.5mm]{0mm}{9mm}
&$\binm{n}{k}\frac{\binm{2n}{k}^2}{\binm{n+k}{k}^2}
\big\{\sst1+(n-2k)(H_k+2H_{n+k})\big\}$
&$(-1)^{\ell}
\binm{n}{\ell}
\frac{\binm{3n+\ell}{\ell}}
     {\binm{n+\ell}{\ell}^2}$
&$\double{\text{Thm\:\ref{thm:2/2}:}
	\:b=0}{c=e=1\:\&\:d\to\infty}$\\
\hline8\rule[-3.5mm]{0mm}{9mm}
&$\binm{n}{k}^4
\frac{\binm{n+k}{k}}{\binm{2n}{k}}
\big\{\sst1+(n-2k)(4H_k-H_{n+k})\big\}$
&$(-1)^{\ell}
\frac{\binm{n}{\ell}\binm{n+\ell}{\ell}^2}
     {\binm{2n}{n}}$
&$\double{\text{Thm\:\ref{thm:2/2}:}
	\:b=1}{c=e=0\:\&\:d\to\infty}$\\
\hline9\rule[-3.5mm]{0mm}{9mm}
&$\binm{n}{k}^3
\frac{\binm{n+k}{k}^2}{\binm{2n}{k}^2}
\big\{\sst1+(n-2k)(3H_k-2H_{n+k})\big\}$
&$\frac{\binm{1+3n}{n}}{\binm{2n}{n}^2}
\times(-1)^{\ell}
\frac{\binm{n}{\ell}\binm{n+\ell}{\ell}^2}
     {\binm{1+2n+\ell}{\ell}}$
&$\double{\text{Thm\:\ref{thm:2/2}:}
	\:b=d=1}{c=0\:\&\:e\to\infty}$\\
\hline10\rule[-3.5mm]{0mm}{9mm}
&$\binm{n}{k}^4
\frac{\binm{n+k}{k}^2}{\binm{2n}{k}^2}
\big\{\sst1+2(n-2k)(2H_k-H_{n+k})\big\}$
&$\frac{\binm{1+3n}{n}}{\binm{2n}{n}^2}
\times(-1)^{\ell}
\frac{\binm{n}{\ell}\binm{n+\ell}{\ell}^3}
     {\binm{1+2n+\ell}{\ell}}$
&$\double{\text{Thm\:\ref{thm:2/2}:}
	\:b=d=1}{c=e=0}$\\
\hline11\rule[-3.5mm]{0mm}{9mm}
&$\binm{n}{k}^4\frac{\binm{2n}{k}}{\binm{n+k}{k}}
\big\{\sst1+(n-2k)(4H_k+H_{n+k})\big\}$
&$(-1)^n\times
\binm{n}{\ell}^2
\frac{\binm{2n+\ell}{\ell}}
     {\binm{n+\ell}{\ell}}$
&$\double{\text{Thm\:\ref{thm:1/3}:}
	\:b\to\infty}{d=1\:\&\:c=e=0}$\\
\hline12\rule[-3.5mm]{0mm}{9mm}
&$\binm{n}{k}^5\frac{\binm{n+k}{k}}{\binm{2n}{k}}
\big\{\sst1+(n-2k)(5H_k-H_{n+k})\big\}$
&$\frac{(-1)^n}{\binm{2n}{n}}\times
\binm{n}{\ell}^2\binm{n+\ell}{\ell}^2$
&$\double{\text{Thm\:\ref{thm:1/3}:}
	\:b=1}{c=d=e=0}$\\
\hline13\rule[-3.5mm]{0mm}{9mm}
&$\binm{n}{k}^3\frac{\binm{2n}{k}^2}{\binm{n+k}{k}^2}
\big\{\sst1+(n-2k)(3H_k+2H_{n+k})\big\}$
&$(-1)^n\times
\frac{\binm{n}{\ell}^2\binm{3n+\ell}{\ell}}
     {\binm{n+\ell}{\ell}^2}$
&$\double{\text{Thm\:\ref{thm:1/3}:}
	\:b\to\infty}{d=0\:\&\:c=e=1}$\\
\hline14\rule[-3.5mm]{0mm}{9mm}
&$\binm{n}{k}
\frac{\binm{2n}{k}^3}{\binm{n+k}{k}^3}
\big\{\sst1+(n-2k)(H_k+3H_{n+k})\big\}$
&$n\times(-1)^\ell
\frac{\binm{n}{\ell}\binm{3n+\ell}{\ell}}
     {(2n-\ell)\binm{n+\ell}{\ell}^2}$
&$\double{\text{Thm\:\ref{thm:1/3}:}
	\:n>0}{b=0\:\&\:c=d=e=1}$\\
\hline15\rule[-3.5mm]{0mm}{9mm}
&$\binm{n}{k}^2\frac{\binm{2n}{k}^3}{\binm{n+k}{k}^3}
\big\{\sst1+(n-2k)(2H_k+3H_{n+k})\big\}$
&$\frac{(-1)^n}{\binm{2n}{n}}
\times\binm{n}{\ell}
\frac{\binm{2n}{\ell}\binm{3n+\ell}{\ell}}
     {\binm{n+\ell}{\ell}^2}$
&$\double{\text{Thm\:\ref{thm:1/3}:}
	\:b\to\infty}{c=d=e=1}$\\
\hline16\rule[-2.5mm]{0mm}{7mm}
&$\binm{n}{k}^5
\big\{\sst1+5(n-2k)H_k\big\}$
&$(-1)^n
\binm{n}{\ell}^2
\binm{n+\ell}{n}$
&$\double{\text{Thm\:\ref{thm:0/4}:}\:b=c=d=0}
	   {e\to\infty;\:\text{cf.\:\citu{paule}{Eq\:5}}}$\\
\hline17\rule[-2.5mm]{0mm}{7mm}
&$\binm{n}{k}^6
\big\{\sst1+6(n-2k)H_k\big\}$
&$(-1)^n
\binm{n}{\ell}^2
\binm{n+\ell}{n}
\binm{2n-\ell}{n}$
&$\double{\text{Thm\:\ref{thm:0/4}:}
	\:}{b=c=d=e=0}$\\
\hline18\rule[-3.5mm]{0mm}{9mm}
&$\binm{n}{k}^5 
\frac{\binm{2n}{k}}{\binm{n+k}{k}}
\big\{\sst1+(n-2k)
(5H_k+H_{n+k})\big\}$
&$(-1)^n\binm{2n}{n}
\times\binm{n}{\ell}^3
\frac{\binm{2n+\ell}{\ell}}
     {\binm{n+\ell}{\ell}
      \binm{2n}{\ell}}$
&$\double{\text{Thm\:\ref{thm:0/4}:}
	\:e=1}{b=c=d=0}$\\
\hline19\rule[-3.5mm]{0mm}{9mm}
&$\binm{n}{k}^4 
\frac{\binm{2n}{k}^2}{\binm{n+k}{k}^2}
\big\{\sst1+2(n-2k)(2H_k+H_{n+k})\big\}$
&$(-1)^n\binm{2n}{n}
\times\frac{\binm{n}{\ell}^3\binm{3n+\ell}{\ell}}
{\binm{n+\ell}{\ell}^2\binm{2n}{\ell}}$
&$\double{\text{Thm\:\ref{thm:0/4}:}
	\:b=d=0}{c=e=1}$\\
\hline20\rule[-3.5mm]{0mm}{9mm}
&$\binm{n}{k}^3
\frac{\binm{2n}{k}^3}{\binm{n+k}{k}^3}
\big\{\sst1+3(n-2k)(H_k+H_{n+k})\big\}$
&$(-1)^n
\frac{\binm{3n}{n}}{\binm{2n}{n}}
\times\binm{n}{\ell}^2
\frac{\binm{2n}{\ell}\binm{3n+\ell}{\ell}}
     {\binm{n+\ell}{\ell}^2\binm{3n}{\ell}}$
&$\double{\text{Thm\:\ref{thm:0/4}:}
	\:b=0}{c=d=e=1}$\\
\hline21\rule[-3.5mm]{0mm}{9mm}
&$\binm{n}{k}^2
\frac{\binm{2n}{k}^4}{\binm{n+k}{k}^4}
\big\{\sst1+2(n-2k)
(H_k+2H_{n+k})\big\}$
&$(-1)^n\frac{\binm{4n}{n}}{\binm{2n}{n}^2}
\times\binm{n}{\ell}\frac{\binm{2n}{\ell}^2
\binm{3n+\ell}{\ell}}
{\binm{4n}{\ell}\binm{n+\ell}{\ell}^2}$
&$\double{\text{Thm\:\ref{thm:0/4}:}
	\:}{b=c=d=e=1}$\\
\hline\hline
\end{tabular}}

\end{document}